\documentclass[10pt,twoside,titlepage]{book}
 \usepackage[cp1251]{inputenc}
\usepackage[russian,ukrainian,english]{babel}
 \usepackage{amsmath}
\usepackage{amsfonts}
\usepackage{latexsym}
\usepackage{amsthm}
\usepackage{amssymb}
\usepackage{isomath}
\usepackage{calc}
\usepackage{indentfirst}
\usepackage{paralist}
\usepackage[numbers,sort&compress]{natbib}

\usepackage{mystyle2}

\newfont{\tenpoint}{cmr10}
\newfont{\eightpoint}{cmr8}
\newfont{\ninepoint}{cmr9}

\makeindex
\newpage
\newtheorem*{proposition*}{Proposition}

\newtheorem*{theorem*}{Theorem}

\newtheorem*{definition*}{Definition}

\newtheorem*{lemma*}{Lemma}

\newtheorem*{corollary*}{Corollary}

\numberwithin{theorem}{section} \numberwithin{proposition}{section}
\numberwithin{definition}{section} \numberwithin{lemma}{section}
\numberwithin{equation}{section} \numberwithin{corollary}{section}
\numberwithin{remark}{section}
\numberwithin{example}{section}

\newcommand{\sep}{/\kern-2pt/ }

 \renewcommand{\thesection}{\arabic{section}}

\allowdisplaybreaks

\usepackage{afterpage}
\emergencystretch=\maxdimen
\setcitestyle{open={([},close={])}}
\usepackage[dotinlabels]{titletoc}
\usepackage{fancyhdr}
\usepackage{titlesec}
\begin{document}

\thispagestyle{empty}

\begin{center}
IVAN FRANKO NATIONAL UNIVERSITY OF LVIV
\end{center}
\begin{center}
IVANO-FRANKIVSK NATIONAL TECHNICAL UNIVERSITY OF OIL AND GAS
\end{center}

\vskip30pt
\begin{flushright}
Monograph
\end{flushright}

\vskip20pt
\begin{center}
Andriy Bandura and Oleh Skaskiv
\end{center}

\begin{flushright}
UDC 517.555
\end{flushright}

\vskip 20pt

\begin{center}
	Entire functions of several variables of bounded index and PDE's
\end{center}

\vskip 100pt

{\footnotesize
\begin{flushright}
\end{flushright}
}

\vskip 110pt

\begin{center}
Lviv -- 2016
\end{center}

\newpage

\begin{flushleft}
\LARGE	\textbf{Contents}
\end{flushleft}
\contentsline {section}{List of notation}{5}
\contentsline {section}{Introduction}{7}
\contentsline {section}{\numberline {1}Main properties of $BLID_{\mathbf {b}}$}{9}
\contentsline {subsection}{\numberline {1.1}Main definitions. The simplest examples}{9}
\contentsline {subsection}{\numberline {1.2}Some properties of the functions from $Q^n_{\mathbf {b}}$}{11}
\contentsline {subsection}{\numberline {1.3}$BLID_{\mathbf {b}}$ and entire functions of one variable}{16}
\contentsline {subsection}{\numberline {1.4}Sufficient sets for $BLID_{\mathbf {b}}$}{17}
\contentsline {subsection}{Problems}{20}
\contentsline {section}{\numberline {2}Characterizations of function of bounded $L$-index in direction}{20}
\contentsline {subsection}{\numberline {2.1}Local behaviour of directional derivative}{21}
\contentsline {subsection}{\numberline {2.2}Description of $L$-index in different directions}{26}
\contentsline {subsection}{\numberline {2.3}Maximum modulus of the entire $BLID_{\mathbf {b}}$}{29}
\contentsline {subsection}{\numberline {2.4}Analogue of Hayman's Theorem}{35}
\contentsline {subsection}{\numberline {2.5}Maximum and minimum modulus}{39}
\contentsline {subsection}{\numberline {2.6}Behaviour of the directional logarithmic derivative}{43}
\contentsline {subsection}{\numberline {2.7}Boundedness of value $L$-distribution in direction}{50}
\contentsline {subsection}{Problems}{53}
\contentsline {section}{\numberline {3}Applications of $BLID_{\mathbf {b}}$}{54}
\contentsline {subsection}{\numberline {3.1}Bounded $L$-index in direction of solutions of partial differential equations}{54}
\contentsline {subsection}{\numberline {3.2}Bounded $L$-index in direction in a bounded domain}{62}
\contentsline {subsection}{\numberline {3.3}Growth and bounded $L$-index in direction of entire solutions of partial differential equations}{64}
\contentsline {subsection}{\numberline {3.4}An example of a function of unbounded index in direction}{70}
\contentsline {subsection}{\numberline {3.5}Bounded $L$-index in direction of some composite functions}{73}
\contentsline {subsection}{\numberline {3.6}$l$-index boundedness for canonical products}{76}
\contentsline {subsection}{\numberline {3.7}$BLID_{\mathbf {b}}$ and entire func\discretionary {-}{}{}tions with ''planar'' zeros}{78}
\contentsline {subsection}{\numberline {3.8}Existence theorems for entire functions of boun\discretionary {-}{}{}ded $L$-index in direction}{81}
\contentsline {subsection}{\numberline {3.9}Growth of entire functions of bounded $L$-index in direction}{84}
\contentsline {subsection}{Problems}{92}
\contentsline {section}{\numberline {4}Entire functions of bounded $\mathbf {L}$-index in joint variables}{93}
\contentsline {subsection}{\numberline {4.1}Main definitions. The simplest examples}{93}
\contentsline {subsection}{\numberline {4.2}Behaviour of derivatives of $B\mathbfit {L}IJ$}{95}
\contentsline {subsection}{\numberline {4.3}Local behaviour of $B\mathbfit {L}IJ$}{100}
\contentsline {subsection}{\numberline {4.4}Boundedness of $L$-index in all directions $\mathbf {e}_j$}{105}
\contentsline {section}{Review of literature}{110}
\contentsline {section}{References}{116}
\contentsfinish

\newpage

\pagestyle{fancy}
\addtolength{\headheight}{\baselineskip}
\renewcommand{\sectionmark}[1]{\markboth{\textbf{\thesection. \ #1}}{}}
\renewcommand{\subsectionmark}[1]{\markright{\textit{\thesubsection. \ #1}}}
\lhead[\leftmark]{\rightmark} \chead {} \rhead {}
\renewcommand{\headrulewidth}{0.6pt}

 \section*{List of notation\markboth{List of notation}{}} \addcontentsline{toc}{section}{List of notation}\label{Poznachennya}
Let $\mathbb{R}_+:=(0;+\infty),$ $F$ be an entire function in $\mathbb{C}^n,$ $L:\mathbb{C}^n\to\mathbb{R}_{+}$ be a continuous function,
$z=(z_1,\ldots,z_n)\in\mathbb{C}^n,$ $\mathbf{b}\in\mathbb{C}^n\setminus \{0\}$ be a given direction, $\mathbf{K}=(k_{1},\ldots,k_{n})\in\mathbb{Z}^{n}_{+},$ $R=(r_1,\ldots,r_n)\in\mathbb{R}^n_{+}.$
\begin{itemize}
\item $\mathbf{K!}=k_{1}!k_{2}!\cdots k_{n}!$ 
\item $\overline{\mathbf{b}}=(\overline{b}_1,\overline{b}_2,\ldots,\overline{b}_n)$ be a conjugate vector to $\mathbf{b}\in\mathbb{C}^n.$
\item $\langle \mathbf{a},\mathbf{b}\rangle=\sum\limits_{j=1}^{n}a_j\overline{b}_j,$ $\mathbf{a}^{\mathbf{b}}=a_{1}^{b_{1}}a_{2}^{b_{2}}
\ldots a_{n}^{b_{n}}$ for $\mathbf{a},$ $\mathbf{b}\in\mathbb{C}^{n}.$

\item $|\mathbf{z}|=\sqrt{|z_{1}|^{2}+|z_{2}|^{2}+\ldots
+|z_{n}|^{2}}$ be the euclidean norm of $\mathbf{z}=(z_{1},\ldots,z_{n})\in\mathbb{C}^{n}.$

\item $\|\mathbf{K}\|=k_1+\ldots+k_n$ for $\mathbf{K}=(k_1,\ldots,k_n)\in\mathbb{Z}^n_{+}.$

\item $D^n(z^0,R)=\{z\in\mathbb{C}^n: |z_j-z_j^0|<r_j, j=1,\ldots,n\}$ be an open polydisc.

\item $D^n[z^0,R]=\{z\in\mathbb{C}^n: |z_j-z_j^0|\leq r_j, j=1,\ldots,n\}$ be a closed polydisc.

\item $T^n(z^0,R)=\{z\in\mathbb{C}^n: |z_j-z_j^0|=r_j, j=1,\ldots,n\}$ be the skeleton of the polydisc.

\item $\frac{\partial^{\|K\|} F}{\partial Z^K}=\frac{\partial^{k_1+\ldots+k_n}F}{\partial z_1^{k_1}\ldots\partial z_n^{k_n}}$ be a partial derivative.

\item $\mathbf{0}=(0,\ldots,0)$ be the zero vector.

\item $\mathbf{e}=(1,\ldots,1), $ \ $\mathbf{e}_j=(0,\ldots,0, \underbrace{1}_{j-\mbox{th place}}, 0,\ldots,0).$

\item for $\eta\geq 0,$ $z\in\mathbb{C}^n,$ $t_0\in\mathbb{C}$
 we define
$$\lambda^{\mathbf{b}}_{1}(z,t_{0},\eta)=\inf\bigg\{\frac{L(z+t\mathbf{b})}
{L(z+t_{0}\mathbf{b})} \colon 
|t-t_{0}|\leq\frac{\eta}{L(z+t_{0}\mathbf{b})}\bigg\},$$
$$\lambda^{\mathbf{b}}_{2}(z,t_{0},\eta)=\sup\left\{\frac{L(z+t\mathbf{b})}
{L(z+t_{0}\mathbf{b})} \colon 
|t-t_{0}|\leq\frac{\eta}{L(z+t_{0}\mathbf{b})}\right\}.$$
\item $\lambda^{\mathbf{b}}_{1}(z,\eta)=\inf\{\lambda^{\mathbf{b}}_{1}(z,t_{0},\eta):
t_{0}\in\mathbb{C}\},$  $\lambda^{\mathbf{b}}_{1}(\eta)=\inf\{\lambda^{\mathbf{b}}_{1}(z,\eta):
z\in\mathbb{C}^{n}\}.$
\item $\lambda^{\mathbf{b}}_{2}(z,\eta)=\sup\{\lambda^{\mathbf{b}}_{2}(z,t_{0},\eta):
t_{0}\in\mathbb{C}\},$ $\lambda^{\mathbf{b}}_{2}(\eta)=\sup\{\lambda^{\mathbf{b}}_{2}(z,\eta):
z\in\mathbb{C}^{n}\}.$
\item ${Q}^{n}_{\mathbf{b}}$ be a class of functions $L$, which for all $\eta\geq 0$ satisfy
 a condition $0<\lambda^{\mathbf{b}}_{1}(\eta)\leq\lambda^{\mathbf{b}}_{2}(\eta)<+\infty;$ $Q\equiv Q^1_1.$
\item For a given $z^0\in\mathbb{C}^n$ we denote $g_{z^0}(t):=F(z^{0}+t\mathbf{b}).$
If $g_{z^0}(t)\not= 0$\  for all $t\in\mathbb{C}$,
 then  $G^{\mathbf{b}}_{r}(F,z^{0}):=\emptyset;$
if $g_{z^0}(t)\equiv 0,$ then
$G^{\mathbf{b}}_r(F,z^0):=\{z^0+t\mathbf{b}\colon
t\in\mathbb{C}\},$
 and if $g_{z^0}(t)\not\equiv 0$ and $a_{k}^{0}$
are  zeros of $g_{z^0}(t)$,
then 
\[G^{\mathbf{b}}_{r}(F,z^{0}):=\bigcup_{k}\left\{z^{0}+t\mathbf{b}\colon
|t-a_{k}^{0}|\leq\frac{r}{L(z^{0}+a^{0}_{k}\mathbf{b})} \right\},\
\ r>0.\]
\item $G^{\mathbf{b}}_{r}(F):=\bigcup_{z^{0}\in\mathbb{C}^{n}}
{G^{\mathbf{b}}_{r}(F,z^{0})}.$
\item $n\big(r,z^{0},t_{0},{1}/{F}\big)=\sum_{|a_{k}^{0}-t_{0}|\leq
r}\,1$ be a counting function of the zero sequence $(a_{k}^{0})$ for $F(z^0+t\mathbf{b})\not\equiv 0.$
\item $M(r,F,z)=\max\{|F(z+t\mathbf{b})|: |t|=r\},$ where $t\in\mathbb{C},$ $z\in\mathbb{C}^n.$
\item $L\asymp L^{*}$ means that for some
$\theta_{1},\theta_{2}\in\mathbb{R}_{+},$
$0<\theta_{1}\leq\theta_{2}<+\infty$ and for all $z\in\mathbb{C}^{n}$ the inequalities 
$\theta_{1}L(z)\leq L^{*}(z)\leq\theta_{2}L(z)$ hold.
\item $g_z(t)=F(z+t\mathbf{b})$ and $l_z(t)=L(z+t\mathbf{b}),$ where $z\in\mathbb{C}^n,$ $t\in\mathbb{C}.$
\item $\frac{\partial F(z)}{\partial
\mathbf{b}}=\sum_{j=1}^{n}\frac{\partial F(z)}{\partial
z_{j}}{b_{j}}=\langle \mathbf{grad}\ F,
\overline{\mathbf{b}}\rangle,$ $\frac{\partial^{k}F(z)}{\partial
\mathbf{b}^{k}}=\frac{\partial}{\partial
\mathbf{b}}\Big(\frac{\partial^{k-1}F(z)}{\partial
\mathbf{b}^{k-1}}\Big), k \geq 2.$
\item $T_m(z,\tau)=T_m(z,\tau,F,L,\mathbf{b})=\frac{1}{m!L^m(\tau)}\big|\frac{\partial^m F(z)}{\partial \mathbf{b}^m} \big|,$
 $T_m(z)=T_m(z,z).$
 \item $BLID_{\mathbf{b}}\equiv $ ``a function of bounded $L$-index in the direction $\mathbf{b}$.
 \item $BlI\equiv$ ''a function of bounded $l$-index``.
 \item $B\mathbfit{L}IJ\equiv$ ``a function of bounded $\mathbf{L}$-index in joint variables''.
 \end{itemize}

In addition to the above we introduce some additional notations  in the following sections.

\newpage

\section*{Introduction\markboth{Introduction}{}} \addcontentsline{toc}{section}{Introduction}

\bigskip

\bigskip

\bigskip

In the modern theory of functions of several complex variables, a leading role is played by the 
theory of entire functions.
 Methods of investigation of entire functions of several complex variables
 can be divided into several groups.
One of them is based on those properties
which can be obtained from the properties of entire functions of one
variable, considering this entire function $F$ as entire function in
each variable separately. Other methods are arised in the study of so-called slice function
i.e. entire functions of one variable $g(\tau)=F(a+b\tau),$
$\tau\in\mathbb{C},$ which is the restriction of the entire function $F$ to an arbitrary
complex line $\{z=a+b\tau: \tau\in\mathbb{C} \},$ \
$a, b \in\mathbb{C}^n.$
This approach is fundamental in our monograph.

B. Lepson \cite{lepson} investigated properties of entire solutions of
linear differential equations and introduced a new subclass of entire
functions so-called functions of bounded index.
This term is used for the entire functions $f$ for which  there
exists $N\in\mathbb{Z}_{+}$ such that for all $p\in\mathbb{Z}_{+}$ and all
$z\in\mathbb{C}$
$$\frac{|f^{(p)}(z)|}{p!}\leq\max\left\{\frac{|f^{(k)}(z)|}{k!}:
0\leq k\leq N\right\}.$$
These functions have been used in the theory
value distribution and differential equations (see bibliography in [\citenum{shah}]).
In particular, every entire function is a function of bounded value distribution if and only if its derivative is a function of bounded index \cite{Hayman}, and every entire solution of the differential equation $f^{(n)}(t)+\sum\limits_{j=0}^{n-1}a_jf^{(j)}(t)=0$ is a function
of bounded index \cite{shahproc}.

G. Fricke and S. Shah investigated an index boundedness of entire solutions of differential equations \cite{shakh}. Later S. Shah \cite{shah} and W. Hayman \cite{Hayman} independently   proved that every entire function of bounded index is a function of exponential type
that its growth is not higher than of normal type of first order.
 M. Salmassi generalized this concept for entire functions of two variables \cite{salmassi,indsalmassi}.

To go beyond the class of entire functions of exponential type A. D. Kuzyk and M. M. Sheremeta ([\citenum{vidlindex}], see also [\citenum{sher}]) for a continuous function $l: \mathbb{R}_{+}\to \mathbb{R}_{+}$  introduced a concept of entire
functions of bounded $l$-index, replacing in the previous definition the quantity 
$ \frac{ |f^{(j)}(z)|}{j!}$ by 
$ \frac{ |f^{(j)}(z)|}{j!l^{j}(|z|)}$.

The multidimensional case is more difficult so there is no such extensive bibliography, as in one-dimensional.
Definition of an entire function of bounded index in several variables was proposed by H. Krishna and S. Shah
in their paper \cite{krishna}.

Properties of bivariate functions of bounded index  were studied in the paper of M. Salmassi \cite{indsalmassi}.
 A concept of the entire function of bounded $\mathbf{L}$-index
in joint variables  was introduced  by M.~M. Sheremeta and M.~T. Bordulyak  \cite{bagzmin}.
These authors (G. Krishna, S. Shah, M. Salmassi, M. Bordulyak, M. Sheremeta) implemented the first approach to transfer the concept of an entire function of bounded index and of
bounded $l$-index of one variable to the class of entire functions
of several variables.
In this case instead of derivatives in the definition, the partial derivatives are considered.

In this way, there was proved
a number of analogues of theorems that describe properties of entire 
functions of bounded $\mathbf{L}$-index and criteria
of boundedness of $\mathbf{L}$-index for entire functions of several variables.
And there were obtained sufficient conditions of  $\mathbf{L}$-index boundedness of entire solutions of some systems of linear differential equations \cite{dysmarta}.
But this approach does not allow to obtain analogues of the one-dimensional characterization of function of bounded $\mathbf{l}$-index in terms of behaviour the logarithmic derivative outside zero sets.
In particular, attempts to investigate of  $\mathbf{L}$-index boundedness for
some important classes of entire functions (for example infinite
products with "planar" zeros) were unsuccessful by technical difficulties.

On the other hand,  this approach fits to study, for example, entire functions of the form
 $F(z)=f_{1}(z_{1})f_{2}(z_{2})\cdots f_{n}(z_{n}),$
$F(z)=f(z_{1}+z_{2}+\cdots+z_{n})$  etc.

Accordingly, the problem arises to consider and to explore an entire function in several variables of bounded $L$-index using a second approach.

We are much indebted to Professor I. E. Chyzhykov, Ivan Franko National University of Lviv, Lviv, Ukraine, for his useful comments and 
encouragement  in the preparation of the monograph.
Finally, the authors thank the referee for  a number of minor corrections and helpful suggestions improving the clarity of the book.

 And we gratefully acknowledge to Associated Prof. M. T. Bordulyak, Ivan Franko National University of Lviv, Lviv, Ukraine,  by courteous 
acquaintance with the content of her dissertation.

\end{document}